\documentclass[12pt]{article}

\usepackage[utf8]{inputenc}
\usepackage[T1,T2A]{fontenc}
\usepackage[russian,english]{babel}
\usepackage{lmodern}

\usepackage{amsmath, amssymb}

\usepackage{graphicx}
\usepackage{placeins} 
\usepackage{caption}
\usepackage{geometry}
\usepackage{setspace}
\usepackage{titlesec}
\titleformat{\section}{\large\bfseries}{\thesection}{1em}{}
\titleformat{\subsection}{\normalsize\bfseries}{\thesubsection}{1em}{}

\usepackage{enumitem}
\setlist{noitemsep}

\usepackage[hidelinks]{hyperref}
\hypersetup{
	pdftitle  = {Geometric Constructions through Ordered Sets},
	pdfauthor = {Liudmyla Morozova},
	pdfcreator= {LaTeX}
}

\begin{document}
	\begin{center}
		{\LARGE \textbf{Geometric Constructions through Ordered Sets}}\\[10pt]
		\textbf{Liudmyla Morozova}\\[4pt]
		Independent researcher, Netanya, Israel\\[2pt]
		\small ORCID: \href{https://orcid.org/0009-0003-0741-9293}{0009-0003-0741-9293}\\[16pt]
	\end{center}
	
	\vspace{0.5em}
	
	\section*{Introduction}
	
This work discusses an approach to solving geometric construction problems in which the given figure is regarded as an element of an ordered set. The core idea is to construct a geometric sequence, or chain, of figures that forms such a set. The first and often nontrivial step is to construct an appropriate sequence in which each successive element is determined by its predecessors. This provides a logical justification for the construction process and naturally leads to the required solution.

To illustrate the proposed approach, three geometric construction problems are considered:
	\vspace{0.5em}
	
	
	\begin{itemize}[leftmargin=1.5em]
		\item \textbf{Proportional segments} --- the traditional and alternative constructions are compared to illustrate the proposed approach;
		\item \textbf{A Gothic detail} --- a figure that does not admit a standard
		direct construction. However, by considering the corresponding set and
		generalizing the problem, we obtain a family of solutions and can carry
		out the required construction;
		\item \textbf{Proportional angles} --- the proposed approach makes it possible to extend the notion of proportional quantities to include angles. To this end, a sequence of angles composed of different numbers of equal parts is constructed.
	\end{itemize}
	
Although these three problems differ substantially, they are not considered
in isolation. Together they reveal a common structural principle: in solving
a construction problem, instead of searching for a single figure, one can
turn to a process that generates an ordered family of related figures. The
resulting chain not only contains the required solution but also makes
visible the logical path by which it arises. This transition from an
individual construction to a generative process is the main subject of
the work.
	
	\section{Proportional Segments}
	
	This chapter presents two approaches to dividing a segment into proportional parts: the traditional approach and an alternative one. It will be shown that the apparent simplicity of this problem is deceptive, as it conceals a deeper meaning.
	
	\subsection{Construction using Thales’ theorem}
	
The task is to divide the segment \(OB_5\) into five equal parts.
To this end, an auxiliary ray is drawn from \(O\), and five equal segments
whose common length is chosen arbitrarily are laid off successively on it.
Together they form the segment \(OA_5\). Thales’ theorem is then applied
(Fig.~\ref{fig:1}a). 

However, this method requires drawing a number of parallel lines.
To carry out the construction, an auxiliary figure is introduced outside
the given configuration. Once the required division has been obtained,
the result is transferred to the segment \(OB_5\) by means of similar
triangles. Thus, the logical development of the solution takes place
outside the figure specified in the problem, and the auxiliary construction
does not form a single developing structure with the given segment.
\subsection{Inclusion of a Segment in an Ordered Set}

The task is again to divide the segment \(OB_5\) into five equal parts.
An auxiliary ray \(a\) is drawn from \(O\), and five equal segments are
laid off successively on it. Their endpoints serve as the centers of five
circles passing through \(O\). To avoid overloading the figure, the centers
of the circles are not labeled separately.

\begin{figure}[htbp]
	\centering
	\includegraphics[width=0.8\linewidth]{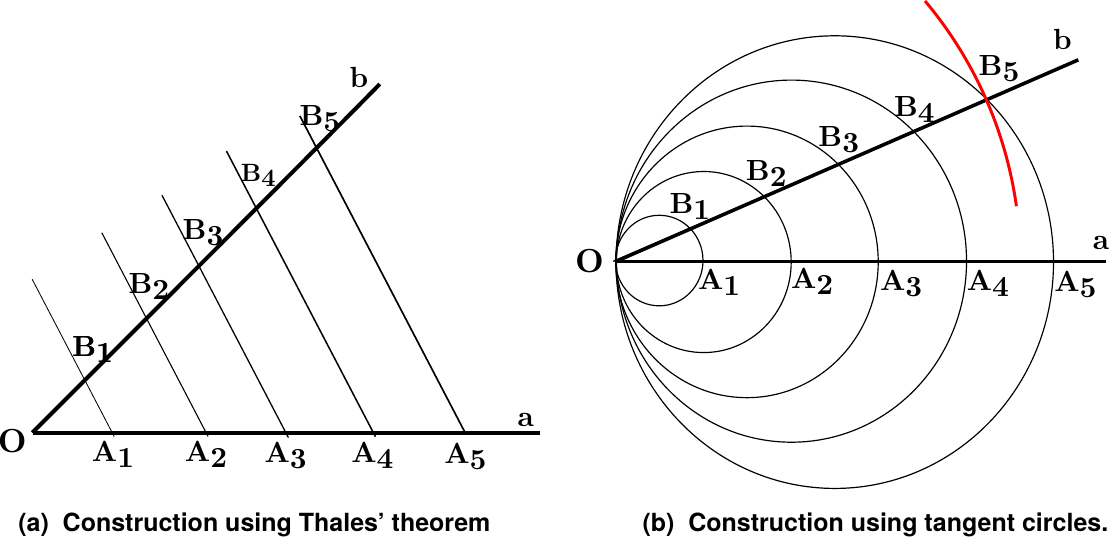}
	\caption{Proportional Segments}
	\label{fig:1}
\end{figure}

The points at which the circles meet the ray \(a\) for the second time
are denoted by \(A_1,A_2,\ldots,A_5\). Thus, the segments
\(OA_1,OA_2,\ldots,OA_5\) are diameters of the constructed circles.
All the circles are internally tangent to one another at \(O\)
(Fig.~\ref{fig:1}b).

An arc centered at \(O\), with radius equal to the length of the given
segment \(OB_5\), is drawn. Its point of intersection with the fifth
circle is denoted by \(B_5\). The segment \(OB_5\) is then drawn. Its
points of intersection with the remaining circles are denoted by
\(B_1,B_2,\ldots,B_4\).

The segments
\[
OA_1,\ OA_2,\ \ldots,\ OA_5
\]
form an external discrete sequence that determines the required
proportions. The constructed circles transfer these relations to the
chords
\[
OB_1,\ OB_2,\ \ldots,\ OB_5.
\]
Since the circles are homothetic with respect to \(O\), the ratios
between the lengths of these chords coincide with the ratios between
the corresponding diameters:
\[
OB_1:OB_2:\ldots:OB_5
=
OA_1:OA_2:\ldots:OA_5
=
1:2:\ldots:5.
\]
Consequently, the points \(B_1,B_2,\ldots,B_4\) divide the given segment
\(OB_5\) into five equal parts.

As the direction of the ray \(OB_5\), on which the constructed chords
lie, varies, their lengths vary continuously, while the ratios between
their lengths remain unchanged. Thus, the segments
\(OB_1,OB_2,\ldots,OB_5\) form a set ordered by length, one element of
which is the given segment \(OB_5\). No parallel lines are required.

The required proportional relation is therefore no longer established
by means of an auxiliary construction outside the given figure; instead,
it arises naturally from the ordered structure itself. Consequently,
the construction possesses both logical completeness and internal
justification.
	\subsection*{Summary}
	
	This chapter presents two approaches to the problem of dividing a segment
	into proportional parts. The first is the traditional construction based
	on Thales’ theorem (Fig.~\ref{fig:1}a). Constructing the required parallel
	lines with a straightedge and compass makes the procedure cumbersome,
	whereas drawing them with a set square takes the solution outside the
	framework of classical geometric construction.
	
	This method uses an auxiliary ray on which five equal segments of an
	arbitrarily chosen length are laid off successively. The resulting division
	is then transferred to the given segment by means of similar triangles.
	However, the auxiliary construction does not arise from the internal
	structure of the original figure and, once the result has been transferred,
	provides no natural continuation. At the same time, Thales’ construction
	already contains an important idea: first, a general scheme of proportional
	relations is created, and then it is applied to the given segment. Here lies
	the seed of the approach developed in the present work.
	
	The second approach uses a chain of tangent circles
	(Fig.~\ref{fig:1}b). It not only eliminates the need to construct
	parallel lines but also brings about a qualitative change in the result.
	The external discrete sequence of diameters \(OA_i\) determines the
	required proportions, and the circles transfer them to the chords
	\(OB_i\). As the direction of the ray \(OB_5\) varies, the lengths of
	these chords vary continuously, while the ratios between their lengths
	remain unchanged. Consequently, instead of dividing a single given
	segment, the construction encompasses a set of segments within the
	prescribed range. In this way, the focus shifts from an individual
	figure to the process of generating an ordered set of figures.
	
\section{A Gothic Detail}

The problem considered here and its algebraic treatment, summarized
below, are taken from Pólya’s \textit{Mathematical Discovery}
[1, pp.~32--34], where the figure is presented as a detail from a
Gothic window. Our geometric approach begins in Section~2.2.

Any problem of geometric construction can be reduced to an algebraic
problem. We shall not present the general theory of such reductions
here; the following example, however, illustrates the idea. A
curvilinear triangle is bounded by the straight line \(AB\) and two
circular arcs, \(AC\) and \(BC\). The center of one circle is at \(A\),
and the center of the other is at \(B\), with each circle passing
through the center of the other. The task is to inscribe in this
curvilinear triangle a circle tangent to all three of its sides;
we shall refer to it as the inner circle.

\subsection{Algebraic Approach}

A configuration of this kind, shown in Fig.~\ref{fig:2}a, is sometimes
found in Gothic architectural ornament. The problem reduces to locating
the center of the inner circle. One locus on which this center must lie
is the axis of symmetry of the region \(ABC\). Another such locus is a
parabola, whose determination requires the methods of analytic geometry.
Thus, this approach is based not on a direct geometric construction,
but on an algebraic determination of the position of the required center.

\subsection{Geometric Approach}

We now complete the arcs of the original figure into full circles and
describe a step-by-step geometric transformation.

\begin{figure}[p]
	\centering
	\includegraphics[width=0.85\linewidth]{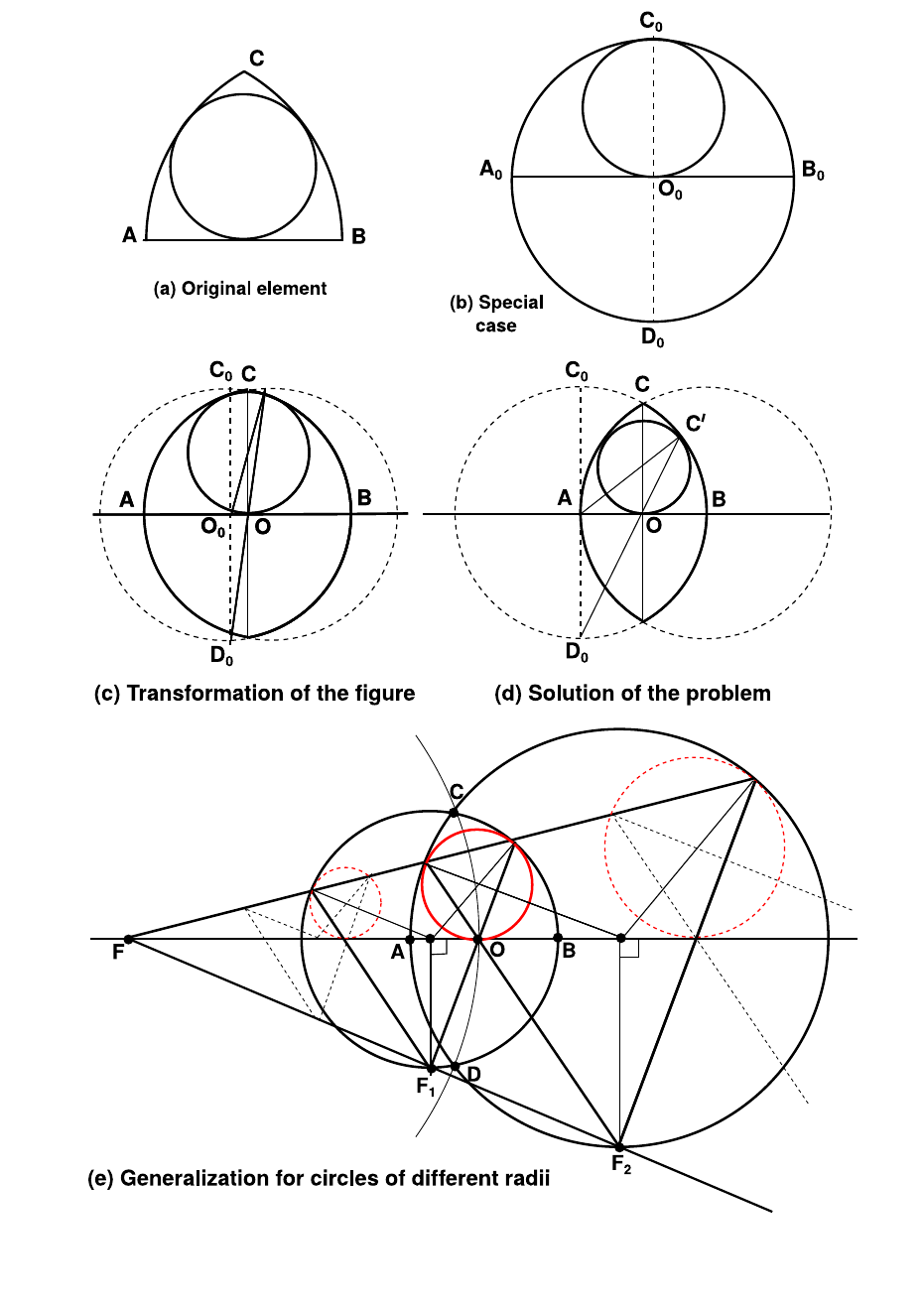}
	\caption{A Gothic Detail}
	\label{fig:2}
\end{figure}

The initial configuration is shown in Fig.~\ref{fig:2}b: the curvilinear
triangle \(A_0B_0C_0\) is formed by a single semicircle. The collinear
points \(D_0\), \(O_0\), and \(C_0\) may be regarded as the vertices
of the degenerate, or ``collapsed,'' triangle \(D_0C_0O_0\).

As the circle is displaced, this triangle begins to unfold. In
Fig.~\ref{fig:2}c, its third vertex is not yet labeled because it lies
too close to the point \(C\), which lies on the chord parallel to the
original diameter \(D_0C_0\). In Fig.~\ref{fig:2}d, this vertex is
denoted by \(C'\), while the point \(A\) occupies the position
corresponding to \(O_0\). Thus, the triangle \(D_0C'A\) represents
the next stage in the development of the initial degenerate triangle
\(D_0C_0O_0\).

The way in which the initial degenerate configuration is interpreted
determines the subsequent construction. Choosing this interpretation
constitutes a nontrivial stage of the approach; once the choice has
been made, the remaining figures arise as a logically ordered sequence.

The displacement produces a second circle, and the curvilinear triangle
\(ABC\) is now bounded by the segment \(AB\) and two arcs, \(AC\) and
\(BC\), each belonging to a different circle. The relation to the point
\(D_0\) is preserved throughout the transformation. In
Fig.~\ref{fig:2}d, the side \(D_0C'\) of the unfolded triangle is also
a diameter of the displaced circle. Its intersection with the segment
\(AB\) determines the point \(O\), the point of tangency of the inner
circle with \(AB\). Thus, the inner circle is constructed not by
locating its center, but by determining its point of tangency with
the segment \(AB\).

The construction is completed in Fig.~\ref{fig:2}d, where the required
solution arises as part of the sequence generated by this transformation.
Since the point \(D_0\) belongs to the initial element of the sequence
but continues to determine the subsequent stages of the construction,
it is, in this sense, a virtual point.

In Fig.~\ref{fig:2}e, the circles move within a sector centered at
the point \(F\). As \(F\) recedes to infinity, the curved trajectories
straighten, and the directions of motion become parallel. Therefore,
the configurations shown in Figs.~\ref{fig:2}c and~\ref{fig:2}d may
be viewed retrospectively as a limiting case of the configuration
shown in Fig.~\ref{fig:2}e. The reverse transition---from parallel
displacement to a sector with a finite center \(F\)---constitutes a
nontrivial step of generalization.

\subsection*{Summary}

Whereas Pólya's algebraic approach reduces the problem to locating
the center of the inner circle, the geometric process developed here
first determines the point \(O\), the point of tangency of this circle
with the segment \(AB\).

The development of the figure begins with the configuration shown in
Fig.~\ref{fig:2}b. The collinear points \(D_0\), \(O_0\), and \(C_0\)
are regarded as the vertices of the degenerate, or ``collapsed,''
triangle \(D_0C_0O_0\). In the initial degenerate configuration, the
point \(C_0\) serves simultaneously as the endpoint of the axis of
symmetry and as a vertex of the collapsed triangle. As the triangle
unfolds, these two geometric roles separate. Figures~\ref{fig:2}c
and~\ref{fig:2}d show two stages of the same process. To avoid
overloading the drawing, the corresponding labels are shown in full
only in Fig.~\ref{fig:2}d.

As the triangle unfolds, the relation to the point \(D_0\) is preserved.
It is precisely the preservation of this relation that determines the
subsequent development of the construction: in Fig.~\ref{fig:2}d,
the diameter \(D_0C'\) intersects the segment \(AB\) at the point \(O\),
thereby determining the point of tangency of the inner circle.

Each panel of Fig.~\ref{fig:2} represents a separate stage in the
development of the problem. Figure~\ref{fig:2}a presents the original
problem. Figure~\ref{fig:2}b identifies the initial configuration from
which the figure can develop; it serves as the beginning of the sequence
corresponding to the given problem. Figure~\ref{fig:2}c shows the process
of development, while in Fig.~\ref{fig:2}d the required construction
arises as a particular result of this process. The position of the
displaced circle is chosen, so the resulting figure represents one of
the possible configurations. Finally, Fig.~\ref{fig:2}e presents a
generalization that follows logically from the original problem.

Thus, the original problem is not exhausted by the construction of a
single solution; it becomes the starting point of a process that
generates an ordered family of interrelated configurations.

\FloatBarrier

\section{Proportional Angles}

This work began with the study of proportional segments. This naturally
raises a related question: can an analogous notion be introduced for
angles? To answer this question, we consider the problem of dividing an
arbitrary angle in a prescribed ratio.

\subsection{A sequence of multiple angles}\label{subsec:sequence-multiples}
We begin with an arbitrarily chosen, sufficiently small initial angle
$\alpha>0$, rather than with the given angle to be divided. The ordered
sequence generated by $\alpha$ will subsequently allow the given angle
to be included in it. Its place in the sequence and its orientation on
the corresponding link will determine the subsequent construction.
The vertex of $\alpha$ is $C_1$, and its rays are highlighted in red
(Fig.~\ref{fig:3}a).

Construct a chain of equal circles of radius $R$ with centers $C_i$
alternating on the two rays of the angle $\alpha$.
By sequentially joining the centers, we obtain a broken line
\[
\langle C_1C_2\cdots C_nC_{n+1}\rangle,
\]
and we obtain the ordered set of multiples
\[
\angle(\alpha),\ \angle(2\alpha),\ \angle(3\alpha),\ \dots,\ \angle(n\alpha).
\]

For the angles supplementary to those between consecutive segments of
the broken line, we obtain:
\begin{align*}
	\angle C_1C_2C_3 &= 180^\circ - 2\alpha, \\
	\angle C_2C_3C_4 &= 180^\circ - 4\alpha, \\
	\angle C_3C_4C_5 &= 180^\circ - 6\alpha, \\
	&\ \ \vdots \\
	\angle C_{n-1}C_nC_{n+1} &= 180^\circ - 2(n-1)\alpha.
\end{align*}

Requiring the angles in the above formulas to be nonnegative yields the basic constraint
\begin{equation*}
	180^\circ - 2(n-1)\alpha \ge 0 \ \Longrightarrow\ 
	\alpha \le \frac{90^\circ}{\,n-1\,}.
\end{equation*}
Thus, as $n$ increases, the maximum admissible value of $\alpha$ decreases.

Instead of dividing the given angle directly, we have constructed an ordered sequence of multiples of $\alpha$. This sequence provides the framework into which the given angle can subsequently be embedded.

\subsection{Doubling the angles of the chain}

Our next task is to learn how to vary the magnitudes of the angles while
preserving the ratios between corresponding angles on different links. We begin
with the simplest special case: we double all the angles. This is
achieved by inscribing within the constructed chain a second chain of
circles whose radii are half those of the original circles (shown in
blue in Fig.~\ref{fig:3}b). The radius can thus be regarded as a
variable parameter.

To avoid overloading the figure, we illustrate this only for the second
and third links (Fig.~\ref{fig:3}b). The numbers (1)--(4) indicate the
successive steps of the construction.

For concreteness, on the third link the angle
$\angle C_4C_3A_4$ is transformed into $\angle C_4C_3D_4$, while on
the second link $\angle C_3C_2A_3$ is transformed into
$\angle C_3C_2E_3$. The corresponding angles preserve their ratio:
\[
\angle C_4C_3D_4 : \angle C_3C_2E_3
=
\angle C_4C_3A_4 : \angle C_3C_2A_3
=
3 : 2.
\]
The red dashed lines indicate the boundaries of possible variation. Thus,
the doubled configuration is not a separate construction but a boundary
case of the structure considered in the next section.

\subsection{Discrete--Continuous Structure}

The first element of the sequence is chosen arbitrarily within the basic
constraint obtained in Section 3.1. Therefore, the same algorithm generates
a family of such sequences.

Their simultaneous presence is reflected in the chosen external chain:
each of its links contains a continuous family of circles corresponding
to different sequences. At the same time, circles of one fixed radius
on all the links form an inner chain.

This gives rise to a two-parameter discrete--continuous structure. The
discrete parameter \(n\) indexes the links of the external chain, while
the continuous parameter \(r\) determines a circle and the corresponding
angle within each link.

The parameter \(r\) varies continuously from \(R\) to \(R/2\). These
boundary values, corresponding to the original and doubled configurations,
were established in Section 3.2. Each intermediate value of \(r\)
corresponds to a definite angle \(\beta\); on the \(n\)-th link, the same
radius corresponds to the angle \(n\beta\).

\begin{figure}[!t]
	\centering
	\includegraphics[width=\linewidth]{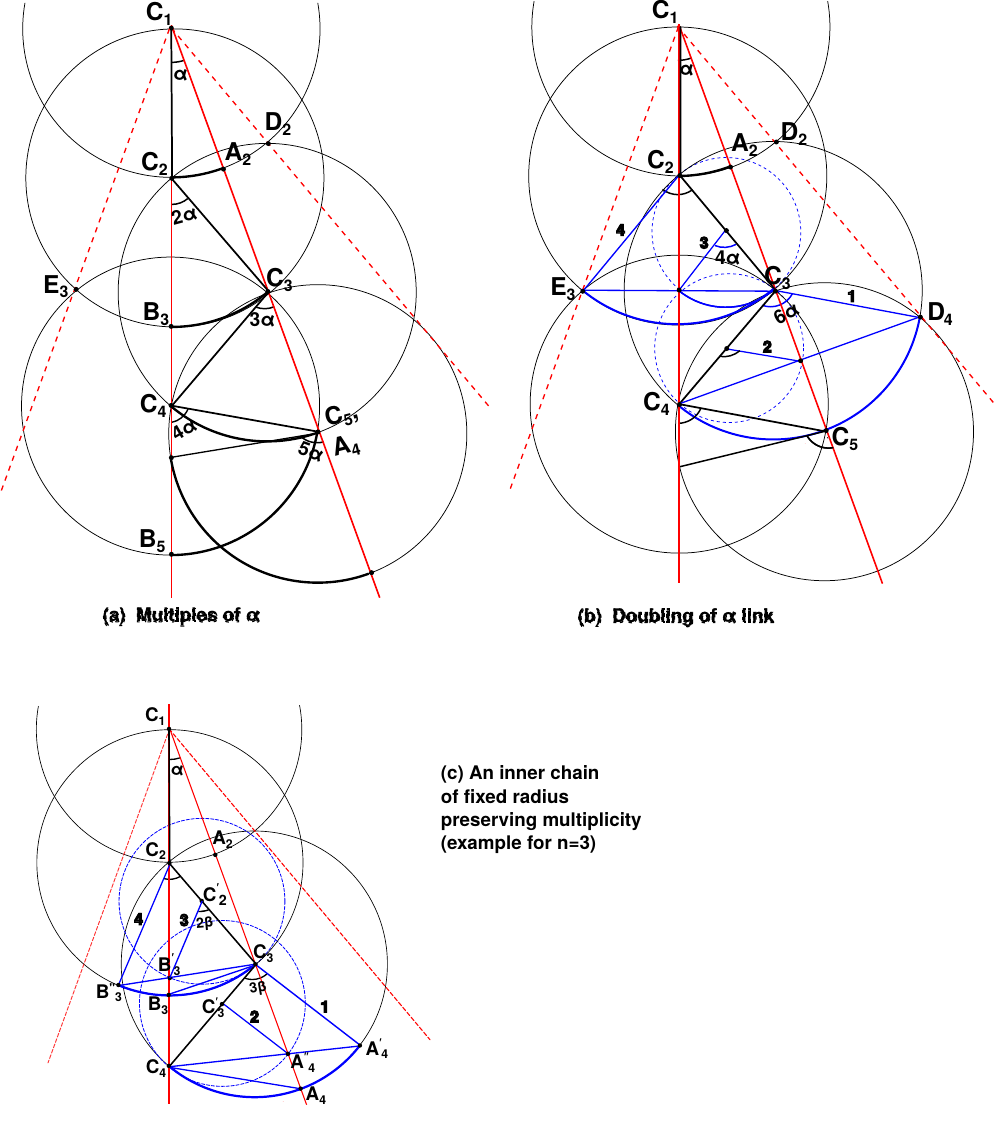}
\captionsetup{justification=raggedright,singlelinecheck=false}
\caption{Proportional Angles. Construction steps in Figs.~3b and 3c:\\
	\hspace*{1.5em}(1) the given angle on the third link;\\
	\hspace*{1.5em}(2) transition to the auxiliary circle;\\
	\hspace*{1.5em}(3) transfer of the construction to the second link;\\
	\hspace*{1.5em}(4) the resulting configuration completing the construction.}
\label{fig:3}
\end{figure}

Depending on the conditions of the problem, we choose an admissible
\(n\)-th link and place the given angle \(n\beta\) on it, as shown in
Fig.~\ref{fig:3}c. The magnitude of this angle belongs to the interval
\((n\alpha,2n\alpha)\). Its position on the chosen link uniquely determines
the required value of the radius \(r_i\), satisfying
\[
\frac{R}{2}<r_i<R.
\]

We thereby select an inner chain passing through all the links and
consisting of circles of radius \(r_i\). In carrying out the construction,
it is sufficient to construct such circles only on the required links.
The ratios between the corresponding angles of this chain remain the
same as in the external chain.

For clarity, consider the case \(n=3\), the simplest example demonstrating
the correspondence between two links. The same procedure extends naturally
as \(n\) increases, although the corresponding drawing soon becomes too
crowded for convenient presentation.

In accordance with the constraint obtained at the end of Section 3.1, for
\(n=3\) we must have \(\alpha\leq45^\circ\). Within this range, choose
\(\alpha\) so that the magnitude of the given angle belongs to the interval
\((3\alpha,6\alpha)\), and place this angle on the third link. Denote its
magnitude by \(3\beta\). Then
\nopagebreak[4]
\[
\alpha<\beta<2\alpha.
\]
Thus, \(\beta\) belongs to the continuous interval associated with the first
link, while the given angle \(3\beta\) selects the corresponding inner chain.

As in Fig.~\ref{fig:3}b, the numbers (1)--(4) in Fig.~\ref{fig:3}c indicate
the successive steps of the construction.

The given angle \(3\beta\) is represented by
\(\angle C_4C_3A'_4\) (Fig.~\ref{fig:3}c). We denote the corresponding
value of the radius \(r_i\) in this case by \(r\). To include the given angle
in the original configuration, construct an inner circle of this radius:
\[
r=A''_4C'_3,\qquad
A''_4C'_3\parallel A'_4C_3.
\]
The inner configuration is thereby inscribed within the outer one.

The radius \(r\), determined on the third link, defines a particular inner
chain. To find the angle corresponding to the prescribed ratio, we pass to
the second link of the same chain.

On the second link, construct a circle of the same radius \(r\), internally
tangent at \(C_3\), with its center lying on the segment \(C_3C'_2\).
Since the radius \(r\) has already been determined on the third link, no
additional parallel line is required. Consequently,
\[
\angle C_4C_3A'_4:\angle C_3C_2B''_3
=
\angle C_4C_3A_4:\angle C_3C_2A_3
=
3:2.
\]
Thus, the angle \(\angle C_3C_2B''_3\) has magnitude \(2\beta\) and
corresponds to the prescribed ratio.

Returning to the original framework, we embed the angle \(2\beta\) into
the second link corresponding to a circle of radius \(R\). This completes
the construction of the required angle and corresponds to step (4) in
Fig.~\ref{fig:3}c.

The same construction applies to every admissible link. At a fixed radius
\(r\), moving between the links of a single inner chain preserves the
proportional relationships between the corresponding angles.

\vspace{-1.5em}
\clearpage
\subsection*{Summary}

In this chapter, we continue the study of proportional relationships
and show that the approach developed in this work applies to angles
as well.

If we begin directly with the given angle, the general problem breaks
down into a series of separate angle-division problems, many of which
cannot be solved with straightedge and compass. Moreover, even when
such a division is possible, individual solutions do not reveal the
overall structure of the problem.

We therefore reverse the direction of the construction. Instead of
beginning with the given angle and seeking one of its parts, we first
choose a sufficiently small angle \(\alpha\) within the range that permits
the construction of the required number of links, and construct an
ordered sequence of its multiples:
\[
\alpha,\ 2\alpha,\ 3\alpha,\ \ldots,\ n\alpha.
\]
Each new element is obtained from its predecessors, and the resulting
sequence forms the external framework for the subsequent solution.

At the next stage, the radius of the circles is treated as a continuously
varying parameter. As the radius decreases from \(R\) to \(R/2\), the
corresponding angle on the \(n\)-th link varies continuously from
\(n\alpha\) to \(2n\alpha\).

For a fixed value of \(\alpha\), the external framework together with the
continuously varying radius forms a two-parameter discrete--continuous
structure:
\begin{itemize}
	\item \(n\) is the discrete parameter indexing the links of the
	external framework;
	\item \(r\) is the continuous parameter determining an inner circle
	and the corresponding value of the angle \(\beta\).
\end{itemize}
\pagebreak
For fixed \(n\), varying \(r\) generates a family of circles of variable
radius in the corresponding link; for fixed \(r\), the corresponding
circles on all the links form an inner chain.

The given angle \(n\beta\) is then placed on an admissible \(n\)-th link.
Its position determines the corresponding value of the radius \(r\) and
thereby selects one of the inner chains. The angles corresponding to its
successive links form the sequence
\[
\beta,\ 2\beta,\ 3\beta,\ \ldots,\ n\beta,
\]
in which the given angle is the \(n\)-th element.

To obtain the required proportional angle, it is sufficient to move to
the required link of the same inner chain and construct there a circle
of radius \(r\). We then return to the original framework and embed the
resulting angle into the corresponding link formed by a circle of radius
\(R\). Thus, the given angle is not divided directly: the solution is
obtained by moving between the links of a single inner chain and then
returning to the original framework. The construction shown in
Fig.~\ref{fig:3}c consists of four steps and requires drawing only one
parallel line.

\clearpage
\section*{Conclusion}
\enlargethispage{2\baselineskip}

The three constructions considered in this work represent successive
stages in the development of the proposed approach. In the first, an
external discrete sequence of segments determines proportional relations,
and the circles transfer these relations to continuously varying chords,
thereby generating an ordered set of figures. In the second, the original
problem becomes the starting point for an ordered family of interrelated
configurations. In the third, a two-parameter discrete--continuous structure arises,
making it possible to extend the concept of proportional magnitudes
to angles.

Thus, embedding a figure into an ordered set not only makes it possible
to trace the logic of its development but also preserves the relations
between successive configurations. The resulting sequence may therefore
be regarded as an information flow.

This flow consists of two complementary components. The first, external,
carries the general idea of the problem, abstracted as far as possible
from its particular conditions---for instance, constructing proportional
magnitudes or inscribing a circle in a curvilinear triangle. Depending
on the problem, it may be represented by a sequence of magnitudes or
by an ordered family of geometric configurations.

The second, internal component is associated with the concrete
realization of this idea. In the first problem, the circles transfer
the relations of the external sequence to continuously varying chords.
In the second, the development of the original configuration leads to
the construction of the corresponding solution. In the third, the given angle determines a particular value of the radius
and thereby selects the corresponding inner chain. Moving to another link
of the same chain makes it possible to determine the required proportional
angle, while returning to the original framework completes its construction.

In the first and third problems, the information flow combines a
discrete sequence with the continuous variation of the associated
magnitudes. In this sense, its structure exhibits a
discrete--continuous duality.
\clearpage
In the third problem, this structure makes it possible to take one
further logical step. Let us return to Fig.~\ref{fig:3}a. The point
$C_2$ may be constructed on either of the two rays of the angle
$\alpha$, and there is no objective reason to prefer one of them.
Here we have freedom of choice. But once the choice has been made,
the broken line $C_1C_2C_3\ldots$ emerges, uniting both rays. This
makes the choice complete.

The red dashed lines indicate the possible positions of the rays of the
angle \(\alpha\) and form a set of three equal angles. This set represents
another manifestation of completeness of choice and sets the limits of
the continuous variation of angles between the original and doubled
configurations shown in Fig.~\ref{fig:3}b.

In turn, each admissible value of the angle \(\alpha\) defines its own
external framework, within which the corresponding family of inner chains
arises. Thus, here too, freedom of choice leads to completeness of choice.
By completeness of choice, we mean the totality of all equally possible
variants.

Combining the notions of information flow and completeness of choice
takes us to the next level of generalization. The relations contained
in the initial configuration are not lost but are redistributed among
the elements of the generated structure. Such redistribution may, in
a limited sense, be regarded as conservation of information within the
ordered structure.
\clearpage
\section*{Limitations}

The approach based on ordered sets does not enlarge the class of
magnitudes constructible with straightedge and compass. It reorganizes
the logical development of a construction within the framework of
classical geometry but does not make it possible to overcome its
inherent limitations.

This boundary becomes particularly clear in the problem of constructing
a segment whose length equals the circumference of a given circle. If
its radius is \(r\), the required segment must have length \(2\pi r\).
Since \(\pi\) is transcendental, such a segment cannot be constructed
with straightedge and compass. The same obstruction arises in the
problem of squaring the circle, whose solution requires constructing
a square with side \(r\sqrt{\pi}\).

Thus, an ordered structure can reveal and organize the relations between
constructible magnitudes, but it cannot make a nonconstructible
magnitude constructible.

\section*{Acknowledgements}
The author gratefully acknowledges the assistance of ChatGPT, a language model
developed by \mbox{OpenAI}, in drafting and refining the English text.

\end{document}